\documentclass[preprint,review,12pt,3p]{elsarticle}

\usepackage{amssymb}
\usepackage[utf8]{inputenc}
\usepackage{multirow}
\usepackage{url}
\usepackage{amsfonts, amsmath, amssymb}
\usepackage{color}
\usepackage{courier}
\usepackage[ruled,vlined]{algorithm2e}
\include{pythonlisting}
\usepackage{lineno}

\emergencystretch=\maxdimen
\newcommand{\norm}[1]{\left\lVert#1\right\rVert}

\journal{Mechanical Systems and Signal Processing}

\begin{document}

\begin{frontmatter}


\title{Robust Partial Quadratic Eigenvalue Assignment Problem: Spectrum Sensitivity Approach}




\author[label1,label2]{José M. Araújo\corref{cor1}}
\address[label1]{On leave from Grupo de Pesquisa em Sinais e Sistemas, Instituto Federal de Educação, Ciência e Tecnologia da Bahia, 41600-270, Salvador, BA, Brazil}

\cortext[cor1]{Corresponding author}
\ead{araujo@ieee.org}

\author[label2]{Carlos E.T. Dórea}
\address[label2]{Universidade Federal do Rio Grande do Norte, Departamento de Engenharia de Computação e Automação, UFRN-CT-DCA, 59078-900 Natal, RN, Brazil}
\ead{cetdorea@dca.ufrn.br}

\author[label2]{Luiz M.G. Gonçalves}
\ead{lmarcos@dca.ufrn.br}

\author[label3]{João B.P. Carvalho}
\address[label3]{ Departamento de Matemática Pura e Aplicada, Universidade Federal do RS, Av Bento Gonçalves 9500, 91509-900, Brazil}

\author[label4]{Biswa N. Datta}
\address[label4]{IEEE Felow; Northern Illinois University, Department of Math. Sciences, 60115, De Kalb, IL, USA}
\ead{dattab@math.niu.edu}

\begin{abstract}
We propose an optimization approach to the solution of the partial quadratic eigenvalue assignment problem (PQEVAP) for active vibration control design with robustness (RPQEVAP). The proposed cost function is based on the concept of sensitivities over the sum and the product of the closed-loop eigenvalues, introduced recently in our paper. Explicit gradient formula for the solutions using state feedback and derivative feedback are derived as functions of a free parameter. These formulas are then used to build algorithms to solve RPQEVAP in a numerically efficient way, with no need to compute new eigenvectors, for both state feedback and  state-derivative feedback designs. Numerical experiments are carried out in order to demonstrate the effectiveness of the algorithms and to compare the proposed method with other methods in the literature, thus showing its effectiveness. 
\end{abstract}

\begin{keyword}
Partial eigenvalue assignment \sep Active vibration control \sep Robustness \sep Spectrum sensitivity
\end{keyword}

\end{frontmatter}


\section{Introduction}
\label{S:1}

Vibrating structures, such as bridges, highways, automobiles, air and space crafts, and others, are usually modeled by a system of second-order differential equations generated by finite element discretization of the original distributed parameter systems. Such second-order system is   known as the Finite Element Model (FEM) in the vibration literature \cite{INMA2001,TIME2001,NANA1993,MING2014,MYRI2006,BETC2013}. These structures sometimes experience dangerous vibrations caused by resonance when excited by external forces including earthquake, gusty winds, weights of human bodies that may result in partial or complete destruction of the structures.
In practice, and very often, such vibrations are controlled by using passive damping forces. Besides being economic to apply it, such an approach has several practical drawbacks: it is ad-hoc in nature and is able to control only localized vibrations. On the other hand, the technique of active vibration control (AVC) is scientifically based and can control vibrations globally in a structure if properly implemented \cite{MORA2006,PALA2016,FULL2008,HERZ2014}. The most important aspect of the AVC implementation is to effectively and efficiently compute the feedback forces needed to control the measured unwanted vibrations, caused by the resonant frequencies.

Recently a mathematically elegant approach that reassigns a few resonant eigenvalues to suitably chosen ones while keeping the other large number of them  and the associated eigenvectors unchanged, has been proposed. This latter approach is known as to having the no-spill over property and the problem of computing the feedback matrices to reassign the unwanted eigenvalues in this way is called PQEVAP (Partial Quadratic Eigenvalue Problem). The approach works exclusively in the second-order setting itself and is capable of  taking advantages of  computationally exploitable inherit structural properties of FEM, such as, definiteness, sparsity, bandness, etc., which are assets in large-scale computational settings.  Typically, the mass and stiffness matrices are symmetric, the mass matrix is positive diagonal and the stiffness matrix is three-diagonal and positive definite or semi-definite.
The most attractive feature of this approach is that the no-spill over property is guaranteed by means of a mathematical theory. This is in sharp contrast with the standard and obvious solution approach of the PQEVAP by transforming a second-order control system to a standard linear state-space. By doing so, one can clearly make use of the existing excellent numerical methods for eigenvalue assignment problems \cite{DATT2003}. However, in this case one needs to deal with a system of dimensions twice that of the original model, which then becomes computationally prohibitive even with a moderate-size model. Notice that the FEM models that arise from practical applications, especially in aerospace and space engineering, and power systems control, could be very large, possibly of multi-million degree  of freedom,  and computational methods for such  large-scale matrix computations are not well-developed \cite{DATT2010}. More importantly,  by transforming to a standard state-space linear system, all the exploitable properties of the FEM, as stated above, will be completely destroyed. By transforming it to a generalized state-space system \cite{DATT2010}, the symmetry can be preserved but not the definiteness. Furthermore, such generalized transformations give rise to descriptor control problems, and  the numerical methods for such control problems, especially for singular and nearly singular and large-scale systems, do not virtually exist \cite{DATT2003}.

A basic solution of the original PQEVAP that meets with the above practical requirements is originally proposed by Datta, et al. \cite{DATT1997}  in  the single input case and then subsequently generalized to the multi-input cases by Datta et al. \cite{DATT2000} and by Elhay and Ram \cite{RAEL2000}. These basic solutions have since then been extended by several authors in recent years using optimization algorithms \cite{BRAH2007,BRAH2009,BAI2010,CAI2010} to compute the two feedback matrices with the important practical properties to ensure that the closed-loop feedback norms are as small as possible and the condition number of the closed-loop eigenvector is minimum.  The associated problems are abbreviated as the MNPQEVAP and RPQPEVAP, respectively. The solutions of MNPQEVAP and   RPQEVAP aim at economic and robust feedback control designs, respectively, both of which are essentials for practical applications. The techniques proposed in those above papers are mostly for state feedback designs which requires explicit knowledge of the state vector. Very often, this vector is not completely available for measurement and, therefore, an observer must be designed \cite{DATT2003}. Unfortunately, numerical algorithms for observer designs for second-order control systems (in second-order setting itself) do not virtually exist \cite{CARV2002}. In order to overcome this difficulty, a few recent papers have been published \cite{ZHAN2015,ZHAN2014,ABDE2013} to solve PQEVAP or full eigenstructure assignment using velocity (state-derivative) feedback. 

In this paper, we propose a new optimization algorithm for RPQEVAP for which the objective function is formulated in terms of the closed-loop mass, stiffness and damping matrices - thus computation of this function and of the associated gradient formula can be performed without explicitly  knowing the closed-loop eigenvectors. This objective function depends upon several spectrum sensitivity results which exhibit these eigenvalue sensitivity relations with the closed-loop feedback matrices \cite{ARAU2016}. The required gradient formulas are derived in the paper in terms of the closed-loop feedback matrices.  These new optimization algorithms are obtained for both cases of the state feedback and state-derivative feedback. Several illustrative examples are given to demonstrate the validity of our results and a comparative study is made with the other methods. 

\section{Preliminary Concepts on Second-Order Systems and the Partial Quadratic Eigenvalue Assignment Problems}

A vibrating structure modeled by a system matrix second-order differential equations has the form:

\begin{equation}
\label{eq1}
\textbf{M}\ddot{\textbf{x}}(t)+\textbf{C}\dot{\textbf{x}}(t)+\textbf{K}\textbf{x}(t)=\mathbf{0},
\end{equation}

\noindent where $\textbf{M},~\textbf{C}$ and $\textbf{K}$ are, respectively, the mass, damping and stiffness matrices, each of order $n$, and $\textbf{x}(t)$ is the displacement vector. Since this model is often generated by using the techniques of finite element,  it is known the Finite Element Model (FEM).
The matrices often have special structures:

\begin{equation}
\textbf{M} =\textbf{ M}^T \succ 0, ~\textbf{C} = \textbf{C}^T,~\textbf{K}=\textbf{K}^T\succeq 0
\end{equation}
The dynamics of such a system are governed by the eigenvalues and eigenvectors of the associated quadratic matrix eigenvalue problem:
\begin{equation}
\label{eq2}
\textbf{Q}(\lambda_k)\textbf{y}_k=0,~k = 1,2...,2n,
\end{equation}

\noindent with the pencil

\begin{equation}
\textbf{Q}(\lambda)=\lambda^2\textbf{M}+\lambda \textbf{C}+\textbf{K} \label{eq3}.
\end{equation}






The details on the quadratic eigenvalue problem can be found in the book \cite{DATT2010}. 
Suppose a control force of the form:
\begin{equation}
\textbf{f}(t)=\textbf{B}\textbf{u}(t)
\end{equation}
\noindent where $\textbf{B}$ is an $n \times m$ control matrix and $\textbf{u}(t)$  is a control vector of $m$  order, applied to the model to control the unwanted vibrations caused by resonances. Thus, we have the control model:

\begin{equation}
\label{eq1c}
\textbf{M}\ddot{\textbf{x}}(t)+\textbf{C}\dot{\textbf{x}}(t)+\textbf{K}\textbf{x}(t)=\mathbf{B}\textbf{u}(t),
\end{equation}Assuming that the state and the velocity vectors $\textbf{x}(t)$ and $\dot{\textbf{x}}(t)$ are known, let's take:

\begin{equation}
\label{eq4}
\textbf{u}(t)=\textbf{F}_s\dot{\textbf{x}}(t)+\textbf{G}_s\textbf{x}(t),
\end{equation}
\noindent where $\textbf{F}_s$ and $\textbf{G}_s$ are two unknown velocity and state feedback matrices. 
Many times the state vector $\textbf{x}(t)$ is not explicitly known, but the velocity vector $\dot{\textbf{x}}(t)$ and the acceleration vector $\ddot{\textbf{x}}(t)$  can rather be estimated. In such case, it is more practical to assume that:
\begin{equation}
\label{eq5}
\textbf{u}(t)=\textbf{F}_d \dot{\textbf{x}}(t)+\textbf{G}_d\ddot{\textbf{x}}(t).
\end{equation}
The control laws defined by (\ref{eq4}) and (\ref{eq5}) , are respectively called the state feedback and the derivative feedback laws. Given then these expressions of the control inputs, the respective closed-loop systems can be written as:

\begin{equation}
\label{eq6}
\textbf{M}\ddot{\textbf{x}}(t) + (\textbf{C}-\textbf{BF}_s)\dot{\textbf{x}}(t)+(\textbf{K}-\textbf{BG}_s)\textbf{x}(t) = 0,
\end{equation}
\begin{equation}
\label{eq7}
(\textbf{M}-\textbf{BG}_d)\ddot{\textbf{x}}(t) + (\textbf{C}-\textbf{BF}_d)\dot{\textbf{x}}(t)+\textbf{K}\textbf{x} (t)= 0.
\end{equation} 
\subsection{Partial Quadratic Eigenvalue Assignment Problem (PQEVAP)}

The partial quadratic eigenvalue assignment problem (PQEVAP) is to assign a few eigenvalues of $Q(\lambda)$, says, $\lambda_1,...,\lambda_p$ ; $p \ll 2n$ , which are considered to cause resonances, to suitably chosen numbers
$\mu_1,...,\mu_p$ by computing the two feedback matrices $\textbf{F}_s$ and $\textbf{G}_s$  for the state feedback case, and $\textbf{F}_d$ and $\textbf{G}_d$ for the derivative feedback case, while leaving the other eigenvalues and eigenvectors unchanged. 


\subsection{Robust  Partial Quadratic Eigenvalue Assignment Problem ( RPQEVAP)}

For practical applications, it is not enough just to compute a pair of feedback matrices. For robust design, one must compute these feedback matrices in such a way that they have norms as small as possible and the closed-loop eigenvalues are as insensitive as possible to small perturbations to data .  
The latter is  equivalent to minimizing the condition number of the closed-loop eigenvector matrix and the problem of finding the feedback matrices such the closed-loop eigenvector matrix has the minimum condition number is called the Robust Partial Quadratic Eigenvalue Assignment Problem or in short, RPQEVAP.
I has been shown recently \cite{BRAH2009,BAI2010} that the eigenvectors matrix $\textbf{Y}_c$ of the closed-loop pencil under partial eigenstructure assignment can be explicitly written as:

\begin{equation}
\mathbf{Y}_c=\begin{bmatrix}
\bar{\textbf{Y}}_1 & \textbf{X}_2\\ 
\bar{\textbf{Y}}_1 \bar{\mathbf{\Lambda}}_1 & \textbf{X}_2 \mathbf{\Lambda}_2
\end{bmatrix}
\end{equation}
\noindent in which the only known quantities are the new assignments  $\bar{\textbf{Y}}_1$ and $\bar{\mathbf{\Lambda}}_1$. Thus, it is a challenge to compute the condition number of the matrix $\textbf{Y}_c$, denoted $\kappa_2(\textbf{Y}_c)$ without having an explicit knowledge of the larger part of this matrix. Some novel ideas have been proposed in the past \cite{BRAH2009,BAI2010} to meet this challenge; these attempts have been made to compute the gradient formulae of the associated optimization problems by knowing only the smaller part of the spectrum and the associated eigenvectors.

In order to state the solutions of these problems in the next section, let's introduce the following notations:

\begin{itemize}
\item $\mathbf{\Lambda}_1 = diag\left (\begin{bmatrix}
\alpha_1 & \beta_1\\ 
-\beta_2 & \alpha_1
\end{bmatrix},\cdots, \begin{bmatrix}
\alpha_l & \beta_l\\ 
-\beta_l & \alpha_l
\end{bmatrix},\lambda_{2l+1},\cdots,\lambda_p\right)$, \\in which $\lambda_k=conj(\lambda_{k+1})=\alpha_k+i\beta_k,~k=1,...,2l$ and $\lambda_{2l+1},...,\lambda_p \in \mathbb{R}$. It is a real representation of the eigenvalues that must to be reassigned. 
\item $\bar{\mathbf{\Lambda}}_1 = diag\left (\begin{bmatrix}
\bar{\alpha}_1 & \bar{\beta}_1\\ 
-\bar{\beta}_2 & \bar{\alpha}_1
\end{bmatrix},\cdots, \begin{bmatrix}
\bar{\alpha}_l & \bar{\beta}_l\\ 
-\bar{\beta}_l & \bar{\alpha}_l
\end{bmatrix},\mu_{2\bar{l}+1},\cdots,\mu_p\right)$, \\in which $\mu_k=conj(\mu_{k+1})=\bar{\alpha}_k+i\bar{\beta}_k,~k=1,...,2\bar{l}$ and $\mu_{2\bar{l}+1},...,\mu_p \in \mathbb{R}$. It is a real representation of the new eigenvalues.
\item $\textbf{Y}_1= \begin{bmatrix}
\Re e(\textbf{y}_1) & \Im m(\textbf{y}_1)  &\cdots   & \Re e(\textbf{y}_l) & \Im m(\textbf{y}_l) & \textbf{y}_{2l+1}\cdots \textbf{y}_p
 \end{bmatrix}$.  It is a real representation of the eigenvectors that must to be reassigned.
 \end{itemize}
Notice that $l$ and $\bar{l}$  are not necessarily  equal, that is, the cardinality  of the complex eigenvalues of the spectrum part to be assigned do not need to equal to that of the reassigned part.

\section{The PQEVAP and RPQEVAP Solutions}

In this section, we first state known solutions to the RPQEVAP and then propose a new optimization approach for the RPQEVAP

\subsection{Solution to PQEVAP}




\begin{itemize}
\item Construction of $\textbf{F}_s$  and $\textbf{G}_s$ : Let be arbitrary $\mathbf{\Gamma}_s \in \mathbb{R}^{m \times p}$ and $\textbf{Z}_s \in \mathbb{R}^{p \times p}$ be the solution of the Sylvester equation:

\begin{equation}
\label{eq9}
\mathbf{\Lambda}_1^T\textbf{Z}_s^T-\textbf{Z}_s^T\bar{\mathbf{\Lambda}}_1=-\textbf{Y}_1^T\textbf{B}\mathbf{\Gamma}_s.
\end{equation}

If $\textbf{Z}_s$ is invertible, and:

\begin{equation}
\label{phis}
\mathbf{\Phi}_s = \mathbf{\Gamma}_s \textbf{Z}_s^{-T},
\end{equation}
\noindent then, it has been shown in \cite{CAI2010} that:

\begin{equation}
\label{sf}
\textbf{F}_s=\mathbf{\Phi}_s \textbf{Y}_1^T\textbf{M},~\textbf{G}_s = \mathbf{\Phi}_s (\mathbf{\Lambda}_1 \textbf{Y}_1^T\textbf{M}+\textbf{C}\textbf{Y}_1^T).
\end{equation}

\item Construction of $\textbf{F}_d$ and $\textbf{G}_d$ : Assume that $0 \notin spec\left(\mathbf{\Lambda}_1\right)$, and let  $\mathbf{\Gamma}_d \in \mathbb{R}^{m \times p}$ and $\textbf{Z}_d$ be the solution of the Sylvester equation:

\begin{equation}
\label{eq10}
\mathbf{\Lambda}_1^T\textbf{Z}_d^T-\textbf{Z}_d^T\bar{\mathbf{\Lambda}}_1=-\mathbf{\Lambda}_1^T\textbf{Y}_1^T\textbf{B}\mathbf{\Gamma}_d.
\end{equation}

If $\textbf{Z}_d$ is invertible, and:
\begin{equation}
\mathbf{\Phi}_d = \mathbf{\Gamma}_d (\textbf{Z}_d^{T}\bar{\mathbf{\Lambda}}_1)^{-1},
\end{equation}

\noindent then it has been shown in\cite{ZHAN2015} that:

\begin{equation}
\label{sd}
\textbf{F}_d = \mathbf{\Phi}_d\mathbf{\Lambda}_1^{T} \textbf{Y}_1^T\textbf{M},~\textbf{G}_d = -\mathbf{\Phi}_d \textbf{Y}_1^T\textbf{K}.
\end{equation}

\end{itemize}

\subsection{RPQEVAP Solution}




\subsubsection{Spectrum sensitivity}
 
In our recent work, the notion of spectrum sensitivity has been introduced \cite{ARAU2016}. More precisely, eight sensitivities related to the perturbations of the sum and the product of the eigenvalues with respect to changes in the system matrices $\textbf{M},~\textbf{C}$, and $\textbf{K}$ were defined as follows:

\begin{equation}
\label{eq11}
\textbf{S}_{\Pi \textbf{K}s}=\frac{\partial \prod{\lambda}^c_s }{\partial \textbf{K}} = \frac{det~(\textbf{K}-\textbf{BG}_s)}{det~\textbf{M}}(\textbf{K}-\textbf{BG}_s)^{-T},
\end{equation}

\begin{equation}
\label{eq12}
\textbf{S}_{\Pi \textbf{M}s}=\frac{\partial \prod{\lambda}^c_s }{\partial \textbf{M}} = -\frac{det~(\textbf{K}-\textbf{BG}_s)}{det~\textbf{M}}\textbf{M}^{-T}, 
\end{equation}

\begin{equation}
\label{eq13}
\textbf{S}_{\Sigma \textbf{C}s}=\frac{\partial \sum{\lambda}^c_s }{\partial \textbf{C}}=-\textbf{M}^{-T},
\end{equation}

\begin{equation}
\label{eq14}
\textbf{S}_{\Sigma \textbf{M}s}=\frac{\partial \sum{\lambda}^c_s }{\partial \textbf{M}}=-\textbf{M}^{-T}(\textbf{C}-\textbf{BF}_s)^T\textbf{M}^{-T},
\end{equation}

\begin{equation}
\label{eq15}
\textbf{S}_{\Pi \textbf{K}d}=\frac{\partial\prod{\lambda}^c_d}{\partial \textbf{K}}=\frac{det~\textbf{K}}{det~(\textbf{M}-\textbf{B}\textbf{G}_d)}\textbf{K}^{-T},
\end{equation}

\begin{equation}
\label{eq16}
\textbf{S}_{\Pi \textbf{M}d}=\frac{\partial\prod{\lambda}^c_d}{\partial \textbf{M}}=-\frac{det~\textbf{K}}{det~(\textbf{M}-\textbf{B}\textbf{G}_d)}(\textbf{M}-\textbf{BG}_d)^{-T}, 
\end{equation}

\begin{equation}
\label{eq17}
\textbf{S}_{\Sigma \textbf{C}d}=\frac{\partial\sum{\lambda}^c_d}{\partial \textbf{C}}=-(\textbf{M}-\textbf{B}\textbf{G}_d)^{-T},
\end{equation}

\begin{equation}
\label{eq18}
\textbf{S}_{\Sigma \textbf{M}d}=\frac{\partial\sum{\lambda}^c_d}{\partial \textbf{M}}=
-(\textbf{M}-\textbf{B}\textbf{G}_d)^{-T}(\textbf{C}-\textbf{B}\textbf{G}_d)^T(\textbf{M}-\textbf{B}\textbf{G}_d)^{-T}.
\end{equation}

In the above formulas on the sensitivities, the subscripts $s$ and $d$ stand, respectively, for the state feedback and the derivative feedback, and $\lambda^c_s$ and $\lambda_d^c$ stand for the closed-loop eigenvalues. in these cases, respectively 


\subsubsection{RPQEVAP with Spectrum Sensitivity}

 Base on the concepts of eigenvalues sensitivities state above, We now formulate the RPQEVAP as follows. For the case of state feedback, we have

Minimize:
\begin{equation}
\label{fs}
 f_s(\mathbf{\Gamma}_s) = \frac{1}{2}w_{1s}\norm{(\textbf{K}-\textbf{BG}_s)^{-T}}_F^2 +\frac{1}{2}w_{2s}\norm{\textbf{M}^{-T}(\textbf{C}-\textbf{BF}_s)^T\textbf{M}^{-T}}_F^2.
\end{equation}

The case of derivative feedback is similar.

Minimize:
\begin{eqnarray}
\label{fd}
f_d(\mathbf{\Gamma}_d) = \frac{1}{2}w_{1d}\norm{(\textbf{M}-\textbf{BG}_d)^{-T}}_F^2\\
\nonumber +\frac{1}{2}w_{2d}\norm{(\textbf{M}-\textbf{BG}_d)^{-T}(\textbf{C}-\textbf{BF}_d)^T(\textbf{M}-\textbf{BG}_d)^{-T}}_F^2.
\end{eqnarray}
Notice that the first term of (\ref{fs}) is related to (\ref{eq11}) and (\ref{eq12}), which concern to sensitivities of the product of the closed-loop eigenvalues with respect to changes in the stiffness and mass matrices. Similarly, the second term of (\ref{fs})  relates the sensitivities of the sum of the closed-loop eigenvalues with respect to damping and mass matrices via (\ref{eq13}) and (\ref{eq21}). Thus, minimization of (\ref{fs}) is related to the minimization of the sensitivities of closed-loop eigenvalues. Similar remarks hold to expression (\ref{fd}). The weights  $w_{1s},~w_{2s},~w_{1d}$ and $w_{2d}$ are  chosen by the designer to avoid dominance in the minimization problems (\ref{fs}) or (\ref{fd}).



In order to minimize (\ref{fs}) and (\ref{fd}),  the corresponding gradient functions must be computed. In the following, we show how to do so. 

\section{Gradient Formulae Construction for State Feedback and derivative Feedback}

\subsection{PQEVAP with State Feedback - SFRPQEVAP}

\textbf{Proposition 1}: \emph{ Suppose that $\mathbf{U},~\mathbf{V}$ are the solutions of the following Sylvester equations:}
 \begin{equation}
 \label{Us}
\bar{\mathbf{\Lambda}}_1 \textbf{U} - \textbf{U} \mathbf{\Lambda}_1^{T} = -\textbf{Z}_s^{-T}\textbf{P}\mathbf{\Upsilon}  \textbf{B}\Gamma \textbf{Z}_s^{-T},
\end{equation}\begin{equation}
\label{V}
\bar{\mathbf{\Lambda}}_1 \textbf{V} - \textbf{V} \mathbf{\Lambda}_1^{T} = -\textbf{Z}_s^{-T}\textbf{Q}\mathbf{\Theta} \textbf{B} \Gamma \textbf{Z}_s^{-T}.
\end{equation}

\emph{and $\mathbf{Z}_s$ is the same as in (\ref{eq9})}

\emph{Then, the gradient $\nabla_{\mathbf{\Gamma}_s}f_{s}$ is given by:}

\begin{eqnarray}
\nabla_{\mathbf{\Gamma}_s} f_s =\left\{\left[\frac{1}{2}\textbf{Z}_s^{-T}\left(\textbf{Q}\mathbf{\Theta}-\textbf{P}\mathbf{\Upsilon}\right)+\frac{1}{2}(-\mathbf{\textbf{V}}+\textbf{U})\textbf{Y}_1^T\right]\textbf{B}\right\}^T,
\end{eqnarray}

\noindent where
\begin{equation}
\mathbf{\Theta} =w_{1s}(\textbf{K}-\textbf{BG}_s)^{-1}(\textbf{K}-\textbf{BG}_s)^{-T} (\textbf{K}-\textbf{BG}_s)^{-1},
\end{equation}
\begin{equation}\mathbf{\Upsilon} = w_{2s}\textbf{M}^{-2}(\textbf{C}-\textbf{BF}_s)^T\textbf{M}^{-2},
\end{equation}
\begin{equation}\textbf{P} = \textbf{Y}_1^T\textbf{M},
\end{equation}
\begin{equation} \label{Q} \textbf{Q} = \mathbf{\Lambda}_1\textbf{Y}_1^{T}\textbf{M}+\textbf{Y}_1^{T}\textbf{C}.
\end{equation}
\textbf{Proof}: The cost function $f(\mathbf{\Gamma}_s)$ given by (\ref{fs}) can be written using matrix traces. Then, from the definition of the gradient of a scalar function of matrices, the differential $\partial f$ must contain some term of the type $tr\left(\nabla_{\Gamma} f^T\partial \mathbf{\Gamma}_s\right)$. By differentiating (\ref{fs}) with respect to $\textbf{F}_s$ and $\textbf{G}_s$ and applying trace properties such as linearity and trace of cyclic permutations for the matrix product, one has:
\begin{eqnarray}
\label{partfs}
 \partial f_s = \frac{1}{2}tr (\mathbf{\Theta} \textbf{B} \partial \textbf{G}_s- \mathbf{\Upsilon} \textbf{B} \partial \textbf{F}_s)+ \frac{1}{2}tr (\textbf{B}^T\mathbf{\Theta}^T \partial \textbf{G}_s^T- \textbf{B}^T\mathbf{\Upsilon}^T \partial \textbf{F}_s^T).
\end{eqnarray}

The first term in (\ref{partfs}) can be expressed as a function of $\partial \mathbf{\Gamma}_s$. Thus by combining (\ref{phis}) and (\ref{sf}) in order to expand $\partial \textbf{F}_s$ and $\partial \textbf{G}_s$, this gives:
\begin{equation}\partial \textbf{F}_s = (\partial \mathbf{\Gamma}_s - \mathbf{\Gamma}_s \textbf{Z}_s^{-T} \partial \textbf{Z}_s^T)\textbf{Z}_s^{-T}\textbf{P},
\end{equation}
\begin{equation}\partial \textbf{G}_s = (\partial \mathbf{\Gamma}_s - \mathbf{\Gamma}_s \textbf{Z}_s^{-T} \partial \textbf{Z}_s^T)\textbf{Z}_s^{-T}\textbf{Q}.
\end{equation}
The differential $\partial \textbf{Z}_s^T$ can be computed by applying a differentiation rule in  (\ref{eq9}):
\begin{equation}
\label{eq20}
\mathbf{\Lambda}_1^T \partial \textbf{Z}_s^T  - \partial \textbf{Z}_s^T \bar{\mathbf{\Lambda}}_1  = -\textbf{Y}_1^T\textbf{B}\partial \mathbf{\Gamma}_s .
\end{equation}
Returning to (\ref{partfs}) and developing the first argument of the trace term leads to:
\begin{eqnarray}
\label{eq21}
\mathbf{\Theta} \textbf{B} \partial \textbf{G}_s- \mathbf{\Upsilon} \textbf{B} \partial \textbf{F}_s =\\
\nonumber \mathbf{\Theta} \textbf{B} \partial \mathbf{\Gamma}_s \textbf{Z}_s^{-T}\textbf{Q}-\mathbf{\Theta} \textbf{B} \mathbf{\Gamma}_s \textbf{Z}_s^{-T} \partial \textbf{Z}_s^T \textbf{Z}_s^{-T}\textbf{Q} -\mathbf{\Upsilon} \textbf{B} \partial \mathbf{\Gamma}_s \textbf{Z}_s^{-T} \textbf{P} + \mathbf{\Upsilon} \textbf{B} \mathbf{\Gamma}_s \textbf{Z}_s^{-T} \partial \textbf{Z}_s^T \textbf{Z}_s^{-T}\textbf{P} .
\end{eqnarray}

Substituting  (\ref{eq20}) into (\ref{partfs}) and using again the properties of the trace function yields:
\begin{eqnarray}
\label{eq22}
\frac{1}{2}tr (\mathbf{\Theta} \textbf{B} \partial \textbf{G}_s- \mathbf{\Upsilon} \textbf{B} \partial \textbf{F}_s) =\\ \nonumber \frac{1}{2}tr[\textbf{Z}_s^{-T}(\textbf{Q} \mathbf{\Theta} - \textbf{P}\mathbf{\Upsilon})\textbf{B}\partial \mathbf{\Gamma}_s]+\frac{1}{2}tr(-\textbf{Z}_s^{-T}\textbf{Q}\mathbf{\Theta} \textbf{B}\mathbf{\Gamma}_s \textbf{Z}_s^{-T}\partial \textbf{Z}_s+\textbf{Z}_s^{-T}\textbf{P}\mathbf{\Upsilon} \textbf{B}\mathbf{\Gamma}_s \textbf{Z}_s^{-T}\partial \textbf{Z}_s^T).
\end{eqnarray}

Next, consider the following solution for the Sylvester equation in   (\ref{eq20}) \cite{HORN1985}:

\begin{equation}
\label{eq40}
\displaystyle
\partial \textbf{Z}_s^T = \sum_{j=0}^{p-1}{\sum_{k=0}^{p-1}{\gamma_{jk}\left(\mathbf{\Lambda}_1^T\right)^j\left(-\textbf{Y}_1^T\textbf{B} \partial \Gamma\right)(\bar{\mathbf{\Lambda}}_1)^k}}.
\end{equation}

\noindent Substituting (\ref{eq40}) into (\ref{eq21}), we obtain:
\begin{eqnarray}
\nonumber \frac{1}{2}tr(-\textbf{Z}_s^{-T}\textbf{Q}\mathbf{\Theta} \textbf{B}\mathbf{\Gamma}_s \textbf{Z}_s^{-T}\partial \textbf{Z}_s^T+\textbf{Z}_s^{-T}\textbf{P}\mathbf{\Upsilon} \textbf{B}\mathbf{\Gamma}_s \textbf{Z}_s^{-T}\partial \textbf{Z}_s^T)=\\
\nonumber \frac{1}{2}tr\left[-\textbf{Z}_s^{-T}\textbf{Q}\mathbf{\Theta} \textbf{B}\mathbf{\Gamma}_s \textbf{Z}_s^{-T}\sum_{j=0}^{p-1}{\sum_{k=0}^{p-1}{\gamma_{jk}\left(\mathbf{\Lambda}_1^T\right)^j\left(-\textbf{Y}_1^T\textbf{B} \partial \mathbf{\Gamma}_s\right)(\bar{\mathbf{\Lambda}}_1)^k}}\right]+\\
\nonumber \frac{1}{2}tr\left[\textbf{Z}_s^{-T}\textbf{P}\mathbf{\Upsilon} \textbf{B}\mathbf{\Gamma}_s \textbf{Z}_s^{-T}\sum_{j=0}^{p-1}{\sum_{k=0}^{p-1}{\gamma_{jk}\left(\mathbf{\Lambda}_1^T\right)^j\left(-\textbf{Y}_1^T\textbf{B} \partial \mathbf{\Gamma}_s\right)(\bar{\mathbf{\Lambda}}_1)^k}}\right]=\\
\nonumber \frac{1}{2}tr\left[-\sum_{j=0}^{p-1}{\sum_{k=0}^{p-1}{\gamma_{jk}(\bar{\mathbf{\Lambda}}_1)^k\left(-\textbf{Z}_s^{-T}\textbf{P}\mathbf{\Theta} \textbf{B}\mathbf{\Gamma}_s \textbf{Z}_s^{-T} \partial \mathbf{\Gamma}_s\right)(\mathbf{\Lambda}_1^T)^j}}\textbf{Y}_1^T\textbf{B}\right]+\\
\nonumber \frac{1}{2}tr\left[\sum_{j=0}^{p-1}{\sum_{k=0}^{p-1}{\gamma_{jk}(\bar{\mathbf{\Lambda}}_1)^k\left(-\textbf{Z}_s^{-T}\textbf{P}\mathbf{\Upsilon} \textbf{B}\mathbf{\Gamma}_s \textbf{Z}_s^{-T} \partial \mathbf{\Gamma}_s\right)(\mathbf{\Lambda}_1^T)^j}}\textbf{Y}_1^T\textbf{B}\right]=\\
\frac{1}{2}tr\left[(-\mathbf{\textbf{V}}+\textbf{U})\textbf{Y}_1^T\textbf{B}\right].
\end{eqnarray}

Thus from (\ref{eq21}), we obtain:
\begin{eqnarray}
\frac{1}{2}tr (\mathbf{\Theta} \textbf{B} \partial \textbf{G}_s- \mathbf{\Upsilon} \textbf{B} \partial \textbf{F}_s)= \frac{1}{2}tr\left\{\left[\textbf{Z}_s^{-T}\left(\textbf{Q}\mathbf{\Theta}-\textbf{P}\mathbf{\Upsilon}\right)+(-\mathbf{\textbf{V}}+\textbf{U})\textbf{Y}_1^T\right]\textbf{B} \partial \mathbf{\Gamma}_s\right\}.
\end{eqnarray}

Since if we have:
\begin{equation}\partial f = tr(\mathbb{J}_1 \partial \textbf{R}+\mathbb{J}_2\partial \textbf{T}),
\end{equation}

\noindent  then $\nabla_\textbf{M} f=\mathbb{J}_1^T$. The proposition is then proved $\blacksquare$.

Based on the gradient formula obtained above, we now state the following algorithm for RPQEVAP with state feedback - SFRPQEVAP.

\begin{algorithm}[H]
\KwIn{~The matrices $\textbf{K},~\textbf{C},~\textbf{M},~\mathbf{\Lambda}_1, \bar{\mathbf{\Lambda}}_1,\textbf{Y}_1$; the maximum number of iterations $maxiter$; the tolerance $\epsilon$}
\KwOut{The feedback matrices $\textbf{F}_s,~\textbf{G}_s$}

\nl \textbf{Step 1:} Set $k=1$ and choose $\mathbf{\Gamma}_s^{(1)} = [\gamma_1...\gamma_p] \in \mathbb{R}^{m \times p}$;\\
\nl \textbf{Step 2:} Compute $\textbf{Z}_s,~\textbf{F}_s,~\textbf{G}_s,$  using (9)-(11);\\
\nl \textbf{Step 3:} Compute $\nabla_{\mathbf{\Gamma}_s} f^{(k)}$ using (\ref{Us})-(\ref{Q});\\
\nl  \textbf{Step 4:} If $\norm{\nabla_{\mathbf{\Gamma}_s} f^{(k)}}_F \leq \epsilon$ or $k = maxiter$, stop. If not, set $k \Leftarrow k+1$ and compute a new $\mathbf{\Gamma}_s^{(k+1)}$  using a gradient-based technique - BFGS, Levenberg-Marquardt or other; return to Step 2.
    \caption{{\bf SFRPQEVAP with Spectrum Sensitivity} \label{Algorithm}}
\end{algorithm}




\vspace{18pt}

\subsection{RPQEVAP with Derivative Feedback - DFRPQEVAP}

\textbf{Proposition 2}:\emph{ Suppose that $\mathbf{U},~\mathbf{V}$  are the solutions of the following Sylvester equations:}
 \begin{equation}
 \label{U}
\bar{\mathbf{\Lambda}}_1 \textbf{U} - \textbf{U} \mathbf{\Lambda}_1^{T} = -\textbf{Z}_d^{-T}\textbf{P}\mathbf{\Upsilon}  B\Gamma \textbf{Z}_d^{-T},
\end{equation}
\begin{equation}
\label{V}
\bar{\mathbf{\Lambda}}_1 \mathbf{V} - \mathbf{V} \mathbf{\Lambda}_1^{T} = -\textbf{Z}_d^{-T}\textbf{Q}\mathbf{\Theta} \textbf{B} \Gamma \textbf{Z}_d^{-T}.
\end{equation}
\noindent \emph{and $\mathbf{Z}_d$ be the same as in (\ref{eq10}).} 

\emph{Then, the gradient $\nabla_{\mathbf{\Gamma}_d}f_{d}$,  (\ref{fd}) is given by:}
\begin{eqnarray}
\nabla_{\mathbf{\Gamma}_d} f_d =\left\{\left[\frac{1}{2}\textbf{Z}_d^{-T}\left(\textbf{Q}\mathbf{\Theta}-\textbf{P}\mathbf{\Upsilon}\right)+\frac{1}{2}(-\mathbf{V}+\textbf{U})\textbf{Y}_1^T\right]\textbf{B}\right\}^T,
\end{eqnarray}

\noindent where
\begin{equation}
\mathbf{\Theta} =w_{1d}(\textbf{M}-\textbf{BG}_d)^{-1}(\textbf{M}-\textbf{BG}_d)^{-T} (\textbf{M}-\textbf{BG}_d)^{-1}+w_{2d} \times
\end{equation}
\begin{equation}
\nonumber [(\textbf{M}-\textbf{BG}_d)^{-1}(\textbf{C}-\textbf{BF}_d)(\textbf{M}-\textbf{BG}_d)^{-1}(\textbf{M}-\textbf{BG}_d)^{-T}
(\textbf{C}-\textbf{BF}_d)^T(\textbf{M}-\textbf{BG}_d)^{-T}(\textbf{M}-\textbf{BG}_d)^{-1}+
\end{equation}
\begin{equation}
\nonumber (\textbf{M}-\textbf{BG}_d)^{-1}(\textbf{M}-\textbf{BG}_d)^{-T}(\textbf{C}-\textbf{BF}_d)^T(\textbf{M}-\textbf{BG}_d)^{-T}(\textbf{M}-\textbf{BG}_d)^{-1}(\textbf{C}-\textbf{BF}_d)(\textbf{M}-\textbf{BG}_d)^{-1}],
\end{equation}
\begin{equation}\mathbf{\Upsilon} = w_{2d}\left[(\textbf{M}-\textbf{BG}_d)^{-1}(\textbf{M}-\textbf{BG}_d)^{-T}(\textbf{C}-\textbf{BF}_d)^T(\textbf{M}-\textbf{BG}_d)^{-T}(\textbf{M}-\textbf{BG}_d)^{-1}\right],
\end{equation}
\begin{equation}\textbf{P} = -\textbf{Y}_1^T\textbf{K},
\end{equation}
\begin{equation}\label{Qd} \textbf{Q} = \mathbf{\Lambda}_1^T\textbf{Y}_1^T\textbf{M}.
\end{equation}

\textbf{Proof}: The proof is similar to that of Proposition 1. It is therefore omitted here.

An algorithm called DFRPQEVAP based on the results of Proposition 2, is now stated bellow for a solution of RPQEVAP in the derivative feedback case.

\begin{algorithm}[ht!]
\KwIn{~The matrices $\textbf{K},~\textbf{C},~\textbf{M},~\mathbf{\Lambda}_1, \bar{\mathbf{\Lambda}}_1,\textbf{Y}_1$; the maximum number of iterations $maxiter$; the tolerance $\epsilon$}
\KwOut{The feedback matrices $\textbf{F}_d,~\textbf{G}_d$}

\nl \textbf{Step 1:} Set $k=1$ and choose $\mathbf{\Gamma}_d^{(1)} = [\gamma_1...\gamma_p] \in \mathbb{R}^{m \times p}$;\\
\nl \textbf{Step 2:} Compute $\textbf{Z}_d,~\textbf{F}_d,~\textbf{G}_d,$  using (12)-(14);\\
\nl \textbf{Step 3:} Compute $\nabla_{\mathbf{\Gamma}_d} f^{(k)}$ using (\ref{U})-(\ref{Qd});\\
\nl  \textbf{Step 4:} If $\norm{\nabla_{\mathbf{\Gamma}_d} f^{(k)}}_F \leq \epsilon$ or $k = maxiter$, stop. If not, set $k \Leftarrow k+1$ and compute a new $\mathbf{\Gamma}_d^{(k+1)}$ using a gradient-based technique as BFGS, Levenberg-Marquardt, or other; return to Step 2;
    \caption{{\bf SDFRPQEVAP with Spectrum Sensitivity} \label{Algorithm}}
\end{algorithm}

\textbf{Comments on Algorithms 1 and 2}: The computational complexity of the proposed algorithms is dominated mainly by the matrix inversions and products necessary to compute the matrices $\mathbf{\Theta}$ and $\mathbf{\Upsilon}$ used in gradient calculation for both techniques. Thus the algorithms are of $O(n^3)$, and therefore efficient. Also, note that both the algorithms can be implemented with the help of only a small number of eigenvalues and eigenvectors that need to be replaced.

\section{Numerical Experiments and Comparisons}

In this section, we present the results of the proposed method with those of other existing methods. Specifically, the following methods are considered for our comparisons:

\begin{enumerate}[I]
\item Proposed state feedback method (Method I - Algorithm I);
\item Proposed derivative feedback (Method II - Algorithm II);
\item The method of Cai et. al \cite{CAI2010} (Method III);
\item The method of Wang \cite{WANG2009} (Method IV);
\item The method of Bai et. al \cite{BAI2010} (Method V).
\end{enumerate}

In the first two experiments,  the matrices $\textbf{M},~\textbf{C}$ and $\textbf{K}$  are perturbed by the following quantities:\[\norm{\Delta \textbf{M}}_F \leq 0.0001 \norm{\textbf{M}}_F,~\norm{\Delta \textbf{C}}_F \leq 0.0001 \norm{\textbf{C}}_F,~ \norm{\Delta \textbf{K}}_F \leq 0.0001 \norm{\textbf{K}}_F.\]

Let $\lambda_j^c$ and $\tilde{\lambda}_j^c$ stand to the jth closed-loop eigenvalue of the unperturbed and the perturbed system, respectively. Then, we define:

\begin{equation}D_{en}=\left[\sum_{j=1}^{2n}\left(\lambda_{j}^c-\tilde{\lambda}_j^c\right)^2\right]^{\frac{1}{2}}.
\end{equation}
\noindent This quantity is the deviation of the perturbed closed-loop eigenvalues from the unperturbed ones. And:


\begin{equation}
R_{ss}=\frac{||\textbf{F}||_\textbf{F}^2+||\textbf{G}||_\textbf{F}^2}{||\textbf{F}_{mn}||_\textbf{F}^2+||\textbf{G}_{mn}||_\textbf{F}^2}.
\end{equation}

\noindent is the relative change to minimum norm feedback.  $||\textbf{F}_{mn}||_\textbf{F}$ and $||\textbf{G}_{mn}||_\textbf{F}$ stand for the minimum feedback norms as computed in the paper of Brahma and Datta \cite{BRAH2007}.

\subsection{Experiment I - Random Example}

In this experiment, we consider a random example of order 5 from MATLAB  \texttt{gallery('randcorr',n)}.

\[\textbf{M}=\begin{bmatrix}
1&0.020074&0.16178&-0.00084629&-0.039004\\
0.020074&1&0.25089&0.090954&0.14549\\
0.16178&0.25089&1&-0.13847&0.0026833\\
-0.00084629&0.090954&-0.13847&1&-0.13832\\
-0.039004&0.14549&0.0026833&-0.13832&1\\
\end{bmatrix},\]
\[\textbf{C}=\begin{bmatrix}
1&-0.044725&-0.093248&-0.16885&0.18645\\
-0.044725&1&0.05047&0.38706&-0.29389\\
-0.093248&0.05047&1&0.0028751&-0.086355\\
-0.16885&0.38706&0.0028751&1&0.034282\\
0.18645&-0.29389&-0.086355&0.034282&1\\
\end{bmatrix},\]
\[\textbf{K} = \begin{bmatrix}
1&-0.63971&-0.16469&0.042341&-0.50555\\
-0.63971&1&0.19923&0.072314&0.49672\\
-0.16469&0.19923&1&0.64109&-0.24001\\
0.042341&0.072314&0.64109&1&-0.403\\
-0.50555&0.49672&-0.24001&-0.403&1\\
\end{bmatrix},\]\[\textbf{B}=\begin{bmatrix}
0.3971&0.9226\\
0.1576&0.4583\\
0.7275&0.7742\\
0.9719&0.3286\\
0.1564&0.3638\\
\end{bmatrix}.\]

The eigenvalues $-0.2551 \pm 1.3772i$ are reassigned to  $-1,~-2$, respectively. The weights $w_{1s}$ and $w_{2s}$ are set to be 1 for Method I.

The results in Table \ref{tab:1} show that $\kappa_2(\textbf{Y}_c)$ is comparable for all the four methods, while the others significant measures of robustness, namely $f(\Gamma_s),~ D_{en}$ and $R_{ss}$ are much better with the Method I than the others.
\begin{table}[tbp]
\centering
\caption{Comparisons of Methods in a random example.}
\label{tab:1}
\begin{tabular}{c|c|c|c|c}
\hline
Method & $f_s(\mathbf{\Gamma}_s)$ & $\kappa_2(\textbf{Y}_c)$ & $D_{en}$ & $R_{ss}$ \\ \hline
I                    & 43.9483                                                   & 170.3181   & 0.0014                                                                                    & 1.7502                                                     \\
III                   & 100.2846                                                   & 116.4019   & 0.0012                                                                                    & 12.9499                                                     \\
IV                  & 44.7336                                                   & 159.3025   & 0.0013                                                                                    & 2.6302                                                     \\
V                   & 389.0345                                                  & 114.7821   & 0.0012                                                                                    & 77.3673                                                    \\ \hline
\end{tabular}
\end{table}


\begin{table}[tbp]
\centering
\caption{Comparison of methods in an example on oil rig.}
\label{tab:2}
\begin{tabular}{c|c|c|c|c}
\hline
Method & $f_s(\mathbf{\Gamma}_s)$ & $\kappa_2(\textbf{Y}_c)$ & $D_{en}$ & $R_{ss}$ \\ \hline
I                    & 0.2979                                               & 4.4416e+004 & 0.2813                                                                     & 1.1355                                                     \\
III                   & 0.3143                                               & 3.5915e+004 & 0.2906                                                                     & 1.8560                                                     \\
IV                  & 3.6190                                               & 6.9402e+006 & 0.9580                                                                     & 101.8630                                                   \\ 
V                  & 0.3524                                               & 9.0145e+004 & 0.2837                                                                     & 3.6009                                                   \\\hline
\end{tabular}
\end{table}

\subsection{Experiment II - An Example on Oil Rig}
\label{sec:4.2}

For this experiment, the matrices $\textbf{M}, ~\textbf{K} \in \mathbb{R}^{66\times 66}$ are obtained from Harwell–Boeing Collection BCSSTRUC1 \cite{LEWI1982} which relate a \textbf{statically condensed oil rig model}. Moreover, we set $\textbf{C} =  \textbf{I}_{66\times 66}$ and $\textbf{B}^T= \begin{bmatrix} \textbf{I}_{2\times 2} & 0_{62\times 2} & -\textbf{I}_{2\times 2}\end{bmatrix}^T$. The eigenvalues $-4.7067 \pm 5.2347i$, $-5.1680 \pm 4.2682i,~-5.2067 + 4.1522i$, of the model are reassigned to the positions: $-6 \pm i,~-6 \pm 2i,~-6 \pm 3i$. The weights are set as $w_{1s}=1$ and $w_{2s}=1e-008$. As seem from Table 2, for the Method I, the condition number for the matrix of closed-loop eigenvectors is better than Method V, while it's comparable with those of Methods III and IV. The other measures are substantially better for the Method I than the others.


\subsection{Experiment III}

The matrices for this experiment were taken from ~\cite{QIXU2005}:

\[\textbf{M} = \textbf{I}_{4\times 4},~~ \textbf{C} =diag([0.5~~0~~0~~0.5]),\]
\[\textbf{K} = \begin{bmatrix}
5 & -5 & 0 & 0\\
-5 & 10 & -5 & 0\\
0 & -5 & 10 & -5\\
0 & 0 & -5 & 6
\end{bmatrix},~~ \textbf{B} = \begin{bmatrix}
1 & 0\\
0 & 1\\
0 & 0\\
0 & 0 
\end{bmatrix}.\]

For this experiment, the eigenvalues $-0.0385 \pm 4.1362i$ are reassigned to $-1 \pm i$, and the perturbations in the system matrices are set as:

\[\norm{\Delta \textbf{M}}_F \leq 0.01 \norm{\textbf{M}}_F,~\norm{\Delta \textbf{C}}_F \leq 0.01 \norm{\textbf{C}}_F,~ \norm{\Delta \textbf{K}}_F \leq 0.01 \norm{\textbf{K}}_F.\]

\noindent The weights considered in this case are $w_{1s}=w_{2s}=1$ for Method I and for Method II$w_{1d}=w_{2d}=1$.

As seen from the results of Table 3, the condition numbers for the eigenvector closed-loop matrix for each of the Methods I and II are smaller than that of Method III. The quantity $D_{en}$ for the Method I is comparable with that of Method III, while for Method II is much better. Figure 1 shows the distribution of re-assigned eigenvalues in both cases of the feedbacks.

\begin{figure}[t!]
\label{fig:cloud}
\centering
\includegraphics[scale=0.8]{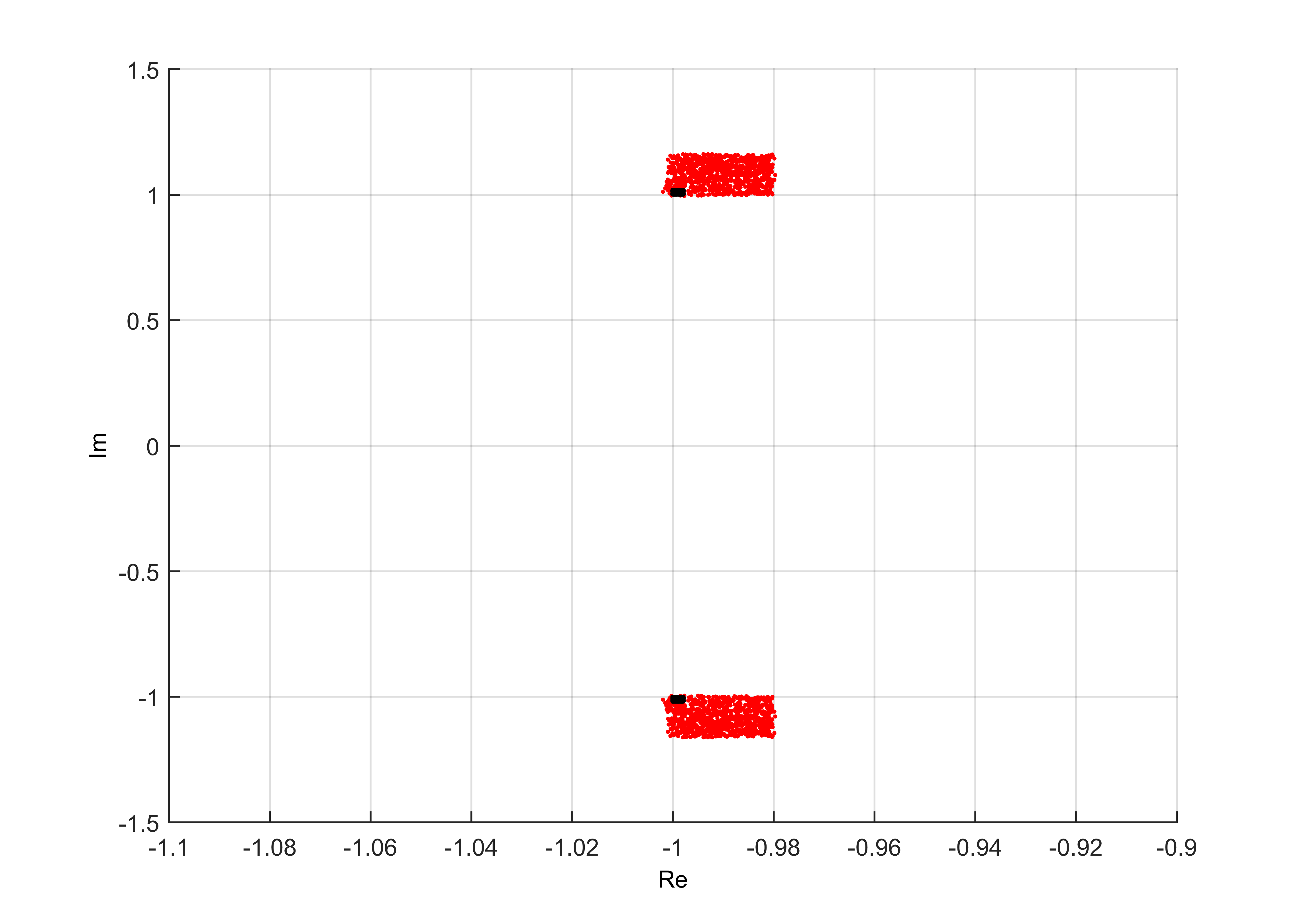}
\caption{Distribution the Reassigned Eigenvalues in the Experiment III with State Feedback (red) and Derivative feedback (black), under Linear Perturbation in Systems Matrices of 1\%.}
\end{figure}

\begin{table}[t!]
\centering
\caption{Condition number and eigenvalue perturbation for the control design in example 3.}
\label{tab:3}
\begin{tabular}{c|c|c}
\hline
Method    & $\kappa_2(\textbf{Y}_c)$ & $D_{en}$ \\ \hline
I   & 21.1073                  & 0.2332                                                                                    \\
II &  46.3772                  &   0.0560                                                                                  \\
III & 72.8761                 & 0.2248                                                                                    \\ \hline
\end{tabular}
\end{table}

\subsection{Experiment IV - Vibration Absorber of a Machine}
In this section, we present the results on system responses of a second-order modeled representing absorber of a machine, taken from \cite{BEAR1996}. 
The matrices $\textbf{M},~\textbf{C}$, and $\textbf{K}$ are given by:

\[\textbf{M} = \textbf{I}_{3\times 3},~~ \textbf{C} = 0\]
\[\textbf{K} = \begin{bmatrix}
2 & 0 & -0.6\\
0 & 2 & -2\\
-0.6 & -2 & 2.68
\end{bmatrix},~~ \textbf{B} = \begin{bmatrix}
1 & 0\\
0 & 0\\
0 & -1 
\end{bmatrix}.\]

The natural natural frequencies of the system are $\pm 2.1108i$, $\pm 1.4142i$, $\pm 0.4737i$.
Now,  a external exciting of the form  $f(t) = 0.1 sin(2.1108t)$ is applied to the system. It's clear that the eigenpair $\pm 2.1108i$ will cause resonance. This eigenpair is then reassigned to $-1 \pm i$  to control the vibration due  to the resonance and the feedback matrices $\textbf{F}_s$ and $\textbf{G}_s$ are computed using Algorithm I. The system responses then are determined for the open-loop system and for the closed-loop system using the Algorithm 1 (Method I) and Method V. These system responses are displayed at Fig. 2. It is seen that the oscillations due to resonance - Fig. 2(a) are well-controlled by applying feedback control forces in both cases - Fig. 2(b),(c).
Finally,  we study the closed-loop system response for Algorithm I and Method V under a perturbation of$+10\%$ and $-10\%$, respectively in matrices \textbf{M} and \textbf{K}. The results are displayed in Fig. 3 . In the Fig. 3(a) the horizontal displacement (red lines) and the torsional tilt (green lines)  under perturbations for the closed-loop system determined by Method V are displayed. The corresponding results are displayed for the closed-loop system obtained by Method I in Fig. 3(b).

\begin{figure}[t!]
\centering
\includegraphics[scale=0.8]{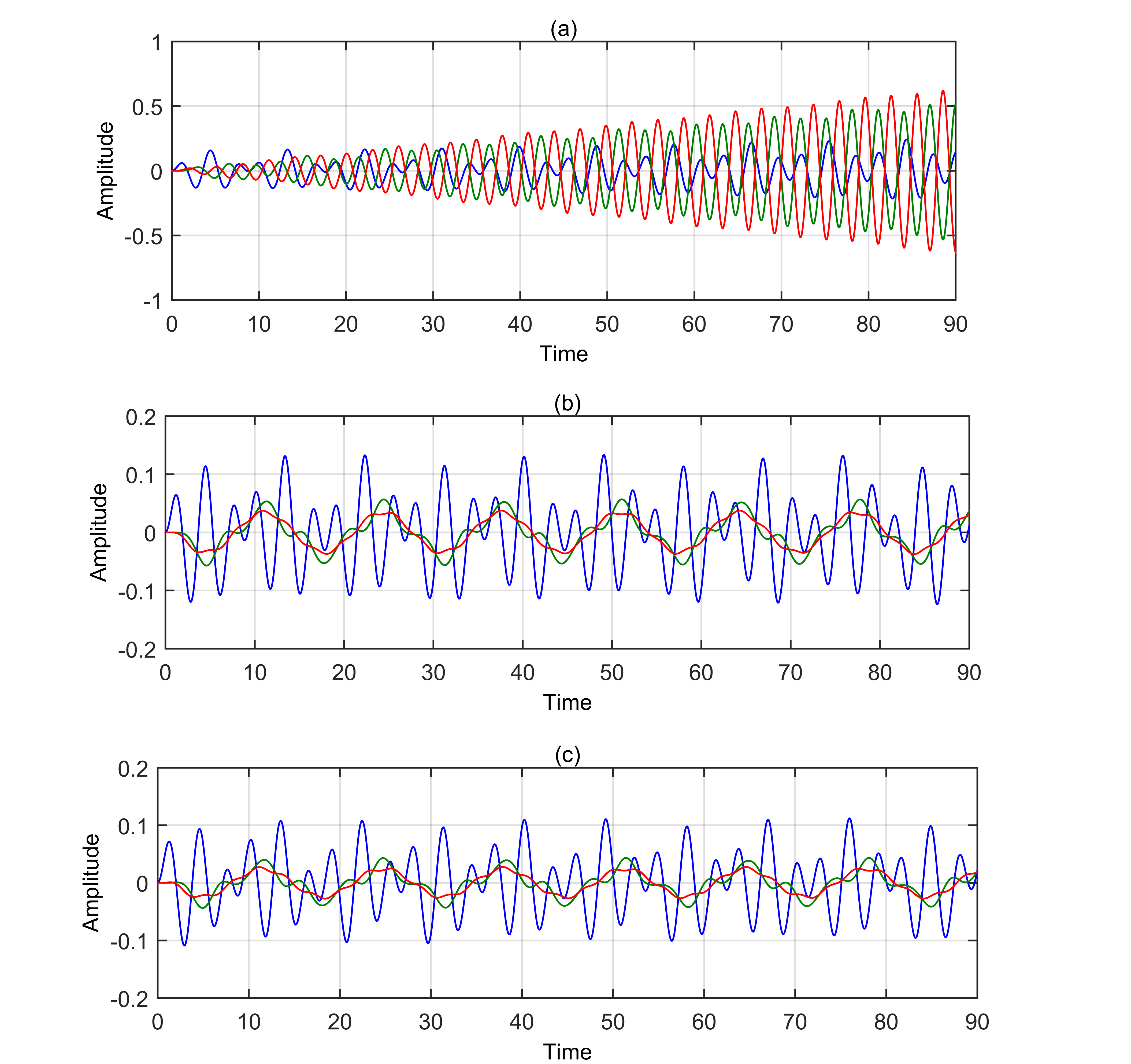} \label{resonance}
\caption{A Study of Controlling Resonant Vibrations by Method I and Method V : (a) Open-loop (b) Closed-loop with Method V (c) Closed-loop with Method I.}
\end{figure}

\begin{figure}[t!]
\centering
\includegraphics[scale=0.8]{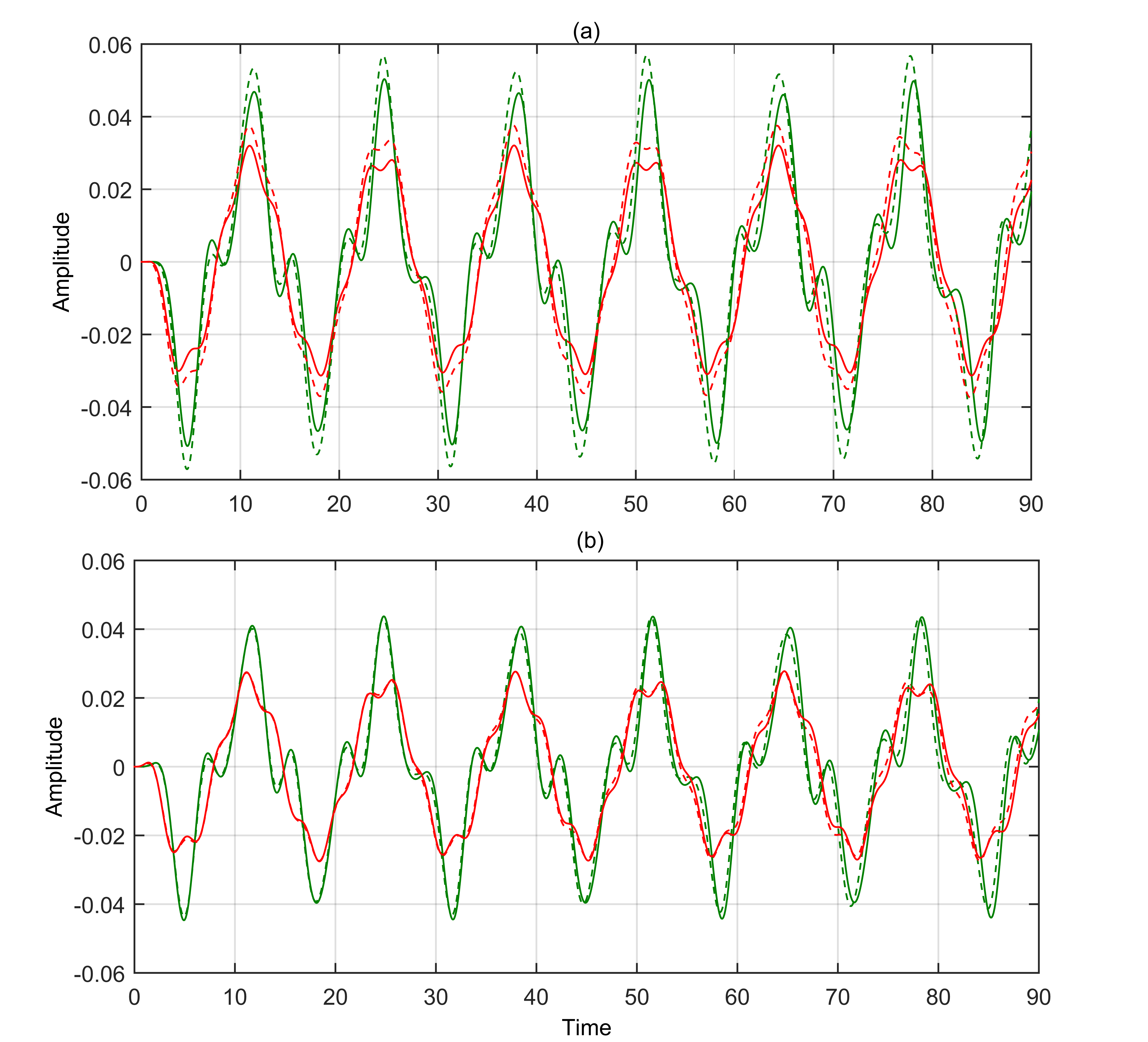} \label{noreson}
\caption{Deviations of the Time Domain Responses for Horizontal Displacement (red lines) and Torsional Tilt (green lines) under Resonant Excitation in Unperturbed  (Continuous) and Perturbed (Dashed) Closed-loop system: (a) Method V (b) Method I.}
\end{figure}

\subsection{Experiment V - Comparison  of the Proposed Algorithms with a Genetic Algorithm}

In this section, we compare the proposed algorithms with a genetic algorithm (GA) which is believed to give a global solution to an optimization problem but heuristic in nature.
The results are displayed in Table 4. Here the superscripts s, d and GA stand for the respective quantities in cases of state feedback, derivative feedback and genetic algorithm. All these three algorithms are applied the three examples in Subsections 5.1, 5.2 and 5.3, considered above.  It is seem that the results on our algorithm are very close or same as those obtained by genetic algorithm for examples 5.1 and 5.3.

\begin{table}[httb!]
\footnotesize
\centering
\caption{Comparison of the proposed gradient-based algorithms 1 and 2 against the meta-heuristic GA optimization.}
\label{tab:IV}
\begin{tabular}{ccccccc}
\hline
    & $f_s(\mathbf{\Gamma}^*)$                                          & $f_d(\mathbf{\Gamma}^*)$                                       & $\norm{\textbf{F}_{I}-\textbf{F}_{GA}}_2$ & $\norm{\textbf{G}_{I}-\textbf{G}_{GA}}_2$ & $\norm{\textbf{F}_{II}-\textbf{F}_{GA}}_2$ & $\norm{\textbf{G}_{II}-\textbf{G}_{GA}}_2$ \\ \hline
5.1 & $\begin{matrix} \mathrm{I} & ~43.9483\\ \mathrm{GA} & 43.9487\end{matrix}$ & -                                                    & 0.0056                    & 0.0410                    & -                         & -                         \\ \hline
5.2 & $\begin{matrix}\mathrm{I} & ~0.2979\\ \mathrm{GA} & 0.2963\end{matrix}$   & -                                                    & 0.6917                    & 3.6533                    & -                         & -                         \\ \hline5.3 & $\begin{matrix}\mathrm{I} & ~16.6393\\ \mathrm{GA} & 16.6451\end{matrix}$ & $\begin{matrix}\mathrm{I} & 2.1972\\ \mathrm{GA} & 2.1972\end{matrix}$ & 0.1062                    & 0.0077                    & 0.0560                    & 0.0097                    \\ \hline
\end{tabular}
\end{table}

\subsection{Experiment VI}

In this experiment, we evaluate the capability of the Algorithm I in reducing the condition number, comparing it with the Method V. We consider the example from \cite{RAEL2000,BAI2010}, with matrices:

\begin{equation}
M =  I_n,~C = 0,~K = \begin{bmatrix}
2 &-1  & 0 & \cdots  & 0 &0 \\ 
-1 & 2 & -1 & \cdots & 0 &0 \\ 
 0& -1 & 2 &\cdots  &0  &0 \\ 
 \vdots & \vdots  & \ddots  & \vdots  & \vdots  &\vdots  \\ 
 0& 0 & \cdots  & -1 & 2 & -1\\ 
 0&  0& \cdots & 0 &  -1& 1
\end{bmatrix},~B=\begin{bmatrix}
I_m\\ 
0
\end{bmatrix}
\end{equation}
\noindent where $n=40$  and $m=3$. The four eigenvalues with smallest absolute value are reassigned to $\lambda_{2k-1}=-k+\sqrt{-10k},~\lambda_{2k}=conj(\lambda_{2k-1})$, $k=1,2$. The weights are chosen as $w_{1s}=0.1$ and $w_{2s}=1$ for Method I.  For sake of comparison, the initial value of $\Gamma_s$  for both methods is taken as:

\begin{equation}
\Gamma_0=\begin{bmatrix}
1 &1  &1  &0 \\ 
 0&1  &1  &1 \\ 
 1&0  &1  &0 
\end{bmatrix}
\end{equation}We then compute the reduction on the condition number for the methods:

\begin{equation}
\Delta \kappa_2(\textbf{Y}_C)\%=100\frac{\kappa^0_{2}-\kappa_2}{\kappa^0_{2}},
\end{equation}
\noindent as well as the quantity $D_{en}$. The results are displayed in Table 4. We observe that, although Method I does not explicitly take into account the condition number in the formulation of the cost function, it gives a reasonable improvement on the condition number after it application,  with a slightly favorable result for the Method V. However, the quantity $D_{en}$ is better in Method I to Method V for perturbations of $1\%$ in both the matrices $\textbf{K}$ and $\textbf{M}$.

\begin{table}[]
\centering
\caption{Reduction for the condition number and deviation of the eigenvalues for Experiment VI. }
\label{my-label}
\begin{tabular}{c|c|c}
\hline
Method & $\Delta \kappa_2(\textbf{Y}_C)\%$ & $D_{en}$ \\ \hline
I      & 49.05\%                           & 0.0412   \\
V      & 63.57\%                           & 0.0451   \\ \hline
\end{tabular}
\end{table}


\section{Conclusions}

A novel approach to RPQEVAP design for second-order controlled linear systems was proposed. The approach consists of minimizing some cost functions that involve sensitivities of the sum and product of closed-loop eigenvalues, named spectrum sensitivities. To this end, new gradient formulae were introduced, and two algorithms were proposed to search for optimal solutions in state feedback and derivative feedback design. In a series of numerical experiments, the proposal was compared with other methods for RPQEVAP solutions, and the results make clear it is competitive to other solutions. Future investigations on this methodology include its application for solving the RPQEVAP in systems with acceleration plus displacement feedback, which is another non-orthodox method for eigenstructure assignment and mode shaping in second-order linear systems.

\section*{Acknowledgments}

The authors would like to thank their Institutions and Brazilian CAPES Foundation for the grants of the research project \#88881.064972/2014-01.

\appendix

\bibliographystyle{elsarticle-num}
\bibliography{lars}







\end{document}